\documentclass[12pt]{amsart}
\usepackage{amsaddr}

\author{ Ming Ho Ng} 
 \address{Department of Mathematics, The Chinese University of Hong Kong\\ Shatin, Hong Kong, P.R. China\\e-mail: mhng@math.cuhk.edu.hk} 
 \author{Yingnan Wang}
 \address{School of Mathematical Sciences, Shenzhen University\\
Shenzhen, Guangdong 518060, P.R. China\\e-mail: ynwang@szu.edu.cn}
\date{\today}
\usepackage{amsmath}
\usepackage{amssymb,latexsym}

\usepackage{color}

\usepackage{enumerate}

\usepackage[T1]{fontenc}

 \usepackage[french,english]{babel}

\makeatletter

\@namedef{subjclassname@2010}{

  \textup{2010} Mathematics Subject Classification}

\makeatother
\newtheorem{thm}{Theorem}[section]

\newtheorem{lem}[thm]{Lemma}

\theoremstyle{definition}

\numberwithin{equation}{section}

\newcommand{\hk}{\mathcal{H}_T}
\newcommand{\sym}{\mathrm{sym}}
\newcommand{\A}{\mathcal{A}}

\newcommand{\ep}{\varepsilon}

\renewcommand{\le}{\leqslant}
\renewcommand{\leq}{\leqslant}
\renewcommand{\ge}{\geqslant}
\renewcommand{\geq}{\geqslant}

\newcommand{\newabstract}[1]{%
  \par\bigskip
  \csname otherlanguage*\endcsname{#1}%
  \csname captions#1\endcsname
  \item[\hskip\labelsep\scshape\abstractname.]
}

\begin{document}

\baselineskip=17pt

\title[]{An omega result for the first sign change of coefficients of symmetric power $L$-functions of Hecke-Maass cusp forms}



\maketitle


\begin{abstract} 
Recently, Lamzouri proved a lower bound for the least positive integer $n_f$ for which the Hecke eigenvalue $\lambda_f(n_f)<0$, showing it satisfies $n_f \ge (\log k)^{1-o(1)}$ for many holomorphic Hecke cusp forms of even integral weight $k$. In this paper, we extend Lamzouri's result to Hecke-Maass forms without assuming the Generalized Ramanujan Conjecture, and further prove that similar results also hold for the coefficients of symmetric power $L$-functions of Hecke-Maass forms.
\end{abstract}

\section{Introduction}
Let $S_k$ be the set of all primitive holomorphic cusp forms of even integral weight $k$ for the full modular group $\mathrm{SL}_2(\mathbb{Z})$. For each $f\in S_k$, let $\lambda_f(n)$ denote the eigenvalue for the $n$th Hecke operator $T_n$. Denote $n_{f}$ by the least positive integer $n$ such that $\lambda_{f}(n)<0$. The study of the first sign change $n_f$ has been a central problem in the theory of holomorphic cusp forms. We refer to \cite{KLSW,Ma12} and the references therein for the history and development of the problem regarding $n_{f}$. In particular, it is well-known that under the generalized Riemann hypothesis (GRH) one can show $n_f\ll (\log k)^2$. By analogy with the least quadratic non-residue, it is folklore to expect that the maximum $\displaystyle \max_{f\in H_k} n_{f}$ should be of order $(\log k)^{1+o(1)}$. Very recently, Lamzouri \cite{La26} obtained an omega result for $n_{f}$ and, moreover, showed that there are at least $\displaystyle{\gg|S_k|\exp\bigg(-3\frac{(\log k)\log\log\log k}{(\log\log k)^2}\bigg)}$ cusp forms $f\in S_k$ satisfying $$n_f\gg \frac{\log k}{\log\log k}.$$

In this short note, we shall consider the analogous problem of the first sign change for Hecke-Maass forms, which carries additional technical difficulties due to the absence of the Generalized Ramanujan Conjecture (GRC). 

Let $\{u_k\}_{k=0}^\infty$ be the set of Hecke-Maass cusp forms for $SL_2(\mathbb{Z})$. Each $u_{k}$ with norm $1$ is an eigenfunction of the Laplacian $\Delta =-y^{2}\left(\frac{\partial ^{2}}{\partial x^{2}}+\frac{\partial ^{2}}{\partial y^{2}}\right)$ with eigenvalue $\frac{1}{4} + t_k^2$ and a simultaneous eigenfunction of all Hecke operators $T_{n}$ with eigenvalue $\lambda_k(n)$. Consider $\hk$ given by $\hk = \{ u_k : t_k \leq T \}$ and we have Weyl's law 
\begin{equation}\label{weyl}
    |\hk| \sim \frac{1}{12} T^2 \quad\textrm{as}\quad T \to \infty.
\end{equation}

Denote by $\alpha_{u_k,1}(p),\alpha_{u_k,2}(p)=\alpha_{u_k,1}^{-1}(p)$ the two Satake parameters of $u_k$ at a prime $p$. The $m$th symmetric power $L$-function associated with $u_k$ is defined for $\Re s>1$ by
\begin{equation}\label{defL}
L(s,\mathrm{sym}^m{u_k})=\prod_p\prod_{\ell=0}^m(1-\alpha_{u_k,1}(p)^{m-2\ell} p^{-s})^{-1}=\sum_{n=1}^\infty\lambda_{\mathrm{sym}^m{u_k}}(n)n^{-s}.
\end{equation}
Denote $n_{\sym^m u_k}$ by the least positive integer $n$ such that $\lambda_{\mathrm{sym}^m{u_k}}(n)<0$. Let $n_{u_k}$ denote $n_{\sym^1 u_k}$. In \cite{LNTW}, $n_{u_k} \ll t_k^{4/5}$ was shown, providing the first explicit exponent for the upper bound. Assuming GRH and GRC, this is expected to be $n_{u_k} \ll (\log t_k)^2$. Furthermore, analogous to the case of holomorphic cusp forms, the magnitude of $\displaystyle \max_{u_k\in \hk} n_{u_{k}}$ is conjectured to be $(\log T)^{1+o(1)}$. In \cite[Corollary 3]{Wang14} Wang showed that $n_{u_k} \ll \log T$ holds for all $u_k \in \hk$, except $u_k$ in an exceptional set of size $\ll |\hk| \exp\left(-c\log T/\log\log T\right)$.

To complement these results, we establish an omega result for $n_{ u_k}$, more precisely for $n_{\sym^m u_k}$.

\begin{thm}\label{Thm:Main}
For sufficiently large $T$, there are $\displaystyle{\gg|\hk|\exp\bigg(-\frac{(\log T)\log\log\log T}{m(\log\log T)^2}\bigg)}$ Hecke-Maass forms $u_k\in \hk$ such that $$n_{\sym^m u_k}\gg \frac{\log T}{m\log\log T}.$$
\end{thm}
This extends the result by Lamzouri \cite{La26} for Hecke eigenvalues of holomorphic Hecke cusp forms to coefficients of symmetric power $L$-functions of Hecke-Maass cusp forms.

\subsection*{Acknowledgements.} 
We would like to thank Professor Yuk-Kam Lau for his comments on the earlier versions of this paper. 

\vskip 4mm
\noindent
{\bf Funding.}
Research of the second author is supported by the National Natural Science Foundation of China (Grant No. 12371006).


\section{Preliminaries}
Gonek and Montgomery \cite[Lemma 7]{GM} introduced the trigonometric
polynomial 
\begin{equation}\label{deff}
    f(\theta)= e^{-i\pi L\theta} \frac{T_L\big(\cos(\pi\theta)/\cos(\pi\delta)\big)}{T_L\big(1/\cos(\pi\delta)\big)},
\end{equation}
where $L$ is a positive integer and
\begin{equation}\label{chebyshevfirst}
T_L(x)= \sum_{m=0\atop m\equiv L\, (\mathrm{mod}\ 2)}^L(-1)^{\frac{L-m}{2}}\, \frac{L}{L+m}\, \binom{\frac{L+m}{2}}{\frac{L-m}{2}}\, 2^{\,m}x^m=\sum_{m=0\atop m\equiv L\, (\mathrm{mod}\ 2)}^L d_m x^m
\end{equation}
is the $L$th Chebyshev polynomial of the first kind. Moreover, they proved 
that for $0<\delta\leq\frac12$ and $\delta\leq\theta\leq1-\delta$,
\begin{equation}\label{boundforf}
|f(\theta)|\leq 2e^{-\pi L\delta}.
\end{equation}
Very recently, Lamzouri \cite[Section 3]{La26} introduced a function
$$
g(\theta)=\bigg|f\Big(\frac{\theta}{2\pi}\Big)\bigg|^2
$$
and proved that for any $\theta\in\mathbb{R}$,
\begin{equation}\label{bound for g1}
    0\le g(\theta)\le 1, \quad
g(0)=1
\end{equation}
and
\begin{equation}\label{Eq:ChebyshevExpansion1}
g(\theta)=\sum_{\ell=0}^{L} a_\ell X_\ell(\theta),
\end{equation}
where 
\begin{equation}\label{Eq:Chebyshevpoly}
X_\ell(\theta)=\frac{\sin(\ell+1)\theta}{\sin\theta}=\frac{e^{i(\ell+1)\theta}-e^{-i(\ell+1)\theta}}{e^{i\theta}-e^{-i\theta}}=\sum_{0\leq j\leq \ell} e^{i(\ell-2j)\theta}
\end{equation}
is the $\ell$th Chebyshev polynomial of the second kind,
$$
a_\ell=\int_0^\pi g(\theta)X_\ell(\theta)\,d\mu_{ST}
=\frac{2}{\pi}\int_0^\pi g(\theta)X_\ell(\theta)(\sin\theta)^2\,d\theta,
$$
and
\begin{align}\label{al}
|a_\ell|
\le
1.
\end{align}

Next, for $\theta\in\mathbb{C}$, we define
\begin{equation}\label{Eq:gdef}
\mathfrak{g}(\theta)=\bigg(e^{i L\theta/2}f\Big(\frac{\theta}{2\pi}\Big)\bigg)^2=\left(\frac{T_L\big(\cos(\theta/2)/\cos(\pi\delta)\big)}{T_L\big(1/\cos(\pi\delta)\big)}\right)^2,
\end{equation}
where $f$ is defined by \eqref{deff}. For $\theta\in\mathbb{R}$, by \eqref{deff}, \eqref{chebyshevfirst} and \eqref{Eq:ChebyshevExpansion1}, we have
$$
\mathfrak{g}(\theta)={g}(\theta)=\sum_{\ell=0}^{L} a_\ell X_\ell(\theta).
$$
Thus, by \eqref{bound for g1}, we have 
\begin{equation}\label{bound for g}
    0\le \mathfrak{g}(\theta)\le 1 \ \textrm{for}\ \theta\in\mathbb{R}\quad\textrm{and}\quad
\mathfrak{g}(0)=1.
\end{equation}
On the other hand, by \eqref{chebyshevfirst} and \eqref{Eq:Chebyshevpoly}, $\mathfrak{g}(\theta)$ and $X_\ell(\theta)$, $\ell\in\mathbb{N}$ are entire functions. Hence, the expansion 
\begin{equation}\label{Eq:ChebyshevExpansion}
\mathfrak{g}(\theta)=\sum_{\ell=0}^{L} a_\ell X_\ell(\theta)
\end{equation} 
extends to all \(\theta \in \mathbb{C}\) by analytic continuation.

We also need the following truncated Kuznetsov trace formula, which is \cite[Lemma 3.1]{LW} (an analogue of \cite[Theorem 6]{BH5}).
\begin{lem}[Kuznetsov]\label{lem3.1}
Let $m,n$ be positive integers. Then we have, for arbitrarily small $\epsilon>0$,
\begin{equation*}
\frac{\pi^2}{T^2}\sum\limits_{t_k\leq T}\alpha_k\lambda_k(m)\lambda_{k}(n)
=\delta_{m,n}+O\left(T^{-1+\epsilon}(mn)^{7/64+\epsilon} + T^{-2}(mn)^{1/4+\epsilon}\right),
\end{equation*}
where $\alpha_k=|\rho_k(1)|^2/\cosh\pi t_k$ and $\delta_{m,n}$ is the Kronecker symbol.
\end{lem}
By \cite[(3.9)]{LW}, we have
\begin{eqnarray*}\label{P}
L(1,{\rm sym}^2u_k) = 2 \alpha_k^{-1}.
\end{eqnarray*}
and by \cite[Main Theorem]{GHL}
$$
L(1,{\rm sym}^2u_k) \gg \frac1{\log (t_k^2+1)}.
$$
Hence, we have
\begin{equation}\label{alphak}
  \alpha_k\ll \log (t_k^2+1).
\end{equation}
\section{Exceptional Hecke eigenvalues}
One main difference between the holomorphic Hecke cusp forms and Hecke-Maass cusp forms is the absence of the Generalized Ramanujan Conjecture (GRC) which predicts that $|\alpha_{u_k,1}(p)|=|\alpha_{u_k,2}(p)|=1$ for any Hecke-Maass cusp form $u_k$ and any prime $p$. The best bound towards GRC
$$
|\alpha_{u_k,j}(p)|\leq p^{\frac7{64}},\quad j=1,2, 
$$ was obtained by Kim and Sarnak \cite{BH14}. Hence, for any Hecke-Maass cusp form $u_k$ and any prime $p$, there exists a unique $\theta_k(p)\in[0,\pi]\cup i(0,\frac7{64}\log p]\cup\pi+i(0,\frac7{64}\log p]$ such that
$$
\lambda_k(p)=2\cos\theta_k(p).$$
It is well-known that
\begin{eqnarray}\label{CP}
\lambda_k(p^n)=X_n(\theta_k(p))=\frac{\sin(n+1)\theta_k(p)}{\sin \theta_k(p)}=\frac{e^{i(n+1)\theta_k(p)}-e^{-i(n+1)\theta_k(p)}}{e^{i\theta_k(p)}-e^{-i\theta_k(p)}}.
\end{eqnarray}

Suppose GRC for $u_k$ fails at some prime $p$. Put
\begin{equation}\label{betadef}
\theta_k(p)=i\beta_k(p)\quad\textrm{or}\quad \theta_k(p)=\pi+i\beta_k(p)\quad\textrm{with}\quad \beta_k(p)\in\mathbb{R}
\end{equation}
according to whether $\theta_k(p)\in i(0,\frac7{64}\log p]$ or $\theta_k(p)\in \pi+i(0,\frac7{64}\log p]$.
Then, by \eqref{Eq:gdef} and \eqref{chebyshevfirst}, we have that for even positive integers $L$,
\begin{align}\label{gpositive}
  \mathfrak{g}(\theta_k(p))&= \frac{1}{T_L\big(1/\cos(\pi\delta)\big)^2} \left(\sum_{m=0\atop m\equiv L\, (\mathrm{mod}\ 2)}^L \frac{d_m}{\big(\cos(\pi\delta)\big)^m} (\cos(\theta_k(p)/2))^m\right)^2\nonumber\\
  &=\frac{1}{T_L\big(1/\cos(\pi\delta)\big)^2} \left(\sum_{m=0\atop m\equiv L\, (\mathrm{mod}\ 2)}^L \frac{d_m}{\big(\cos(\pi\delta)\big)^m} \left(\frac{e^{i\theta_k(p)/2}+e^{-i\theta_k(p)/2}}{2}\right)^m\right)^2\nonumber\\
  &\geq0,
\end{align}
no matter whether $\theta_k(p)\in i(0,\frac7{64}\log p]$ or $\theta_k(p)\in \pi+i(0,\frac7{64}\log p]$.
Moreover, by \eqref{CP}, we have
\begin{align}
|\lambda_k(p^n)|=|X_n(2\cos\theta_k(p))|&=\frac{e^{(n+1)|\beta_k(p)|}-e^{-(n+1)|\beta_k(p)|}}{e^{|\beta_k(p)|}-e^{-|\beta_k(p)|}}\label{eq:1}\\ 
&=
\sum_{0\leq j\leq n} e^{(n-2j)|\beta_k(p)|}.\label{eq:2}
\end{align}

Let $\delta_p= (\log p)/(16m \log \log T)$ and $L_p=4\lfloor m
(\log\log  T)^2/\log p\rfloor+4 \in 4 \mathbb{N}$. We also define $\varepsilon:=4(\log T)^{-\pi/2}$.
Let $u_k\in\hk$, $2\leq z\leq \log T$ be a real number and $H\geq \max\limits_{p\leq z}L_p^{2026}$ be a positive integer to be chosen. Define
\begin{equation}\label{Eq:DefG}
\mathcal{G}(u_k):=\prod_{p\le z} \mathfrak{g}_p(\theta_k(p))
-\sum_{\substack{q\le z\\ q\text{ prime}}}\left(\varepsilon+\frac{\lambda_k(q^H)^2}{H^2}\right)
\prod_{\substack{p\le z\\ p\ne q}} \mathfrak{g}_p(\theta_k(p)).
\end{equation}
Here $\mathfrak{g}_p(\theta_k(p))$ is defined by \eqref{Eq:gdef} with $\delta=\delta_p$ and $L=L_p$. By \eqref{Eq:ChebyshevExpansion}, 
$$
\mathfrak{g}_p(\theta)=\sum_{\ell=0}^{L_p} a_{\ell,p} X_\ell(\theta).
$$
\noindent\textbf{\textrm{Case 1.}} Note that if
$
\theta_k(q)\in [2\pi \delta,\pi]
$
for some prime $q\le z$, then by \eqref{boundforf} and \eqref{Eq:gdef}, we have $0\leq \mathfrak{g}_q(\theta_k(q))\le\varepsilon$ and so combining with \eqref{bound for g} and \eqref{gpositive}, we have 
$$
\mathcal{G}(u_k)
\le
\left(\mathfrak{g}_q(\theta_k(q))-\left(\varepsilon+\frac{\lambda_k(q^H)^2}{H^2}\right)\right)\prod_{\substack{p\le z\\ p\ne q}} \mathfrak{g}_p(\theta_k(p))
\le 0.
$$
\noindent\textbf{\textrm{Case 2.}}
Suppose GRC for $u_k$ fails at some prime $q\leq z$ with $|\beta_k(q)|\geq\frac{L_q+2\log H}{2H-L_q-1}$ (see \eqref{betadef} for the definition of $\beta_k(q)$). Then by \eqref{eq:1}
$$
\frac{\lambda_k(q^H)^2}{H^2}\gg H^{-2}e^{2H|\beta_k(q)|}\max(1,|\beta_k(q)|^{-2})
$$
and by \eqref{eq:1}, \eqref{Eq:ChebyshevExpansion} and \eqref{al},
\begin{align*}
\mathfrak{g}_q(\theta_k(q))  \leq \sum_{\ell=0}^{L_q} |a_{\ell,q}| |X_\ell(\theta_k(q))|  \ll L_q e^{(L_q+1)|\beta_k(q)|}\max(1,|\beta_k(q)|^{-1}).
\end{align*}
Then there exists $T_1>0$ such that for any $T>T_1$ and $|\beta_k(q)|\geq\frac{L_q+2\log H}{2H-L_q-1}$, we have
$$
\mathfrak{g}_q(\theta_k(q))\leq\frac{\lambda_k(q^H)^2}{H^2}
$$
and so combining with \eqref{bound for g} and \eqref{gpositive}, we have 
$$
\mathcal{G}(u_k)
\le
\left(\mathfrak{g}_q(\theta_k(q))-\frac{\lambda_k(q^H)^2}{H^2}\right)\prod_{\substack{p\le z\\ p\ne q}} \mathfrak{g}_p(\theta_k(p))
\le 0.
$$
\noindent\textbf{\textrm{Case 3.}} Suppose GRC for $u_k$ fails at some prime $q\leq z$ with $|\beta_k(q)|<\frac{L_q+2\log H}{2H-L_q-1}$. If $\theta_k(q)\in i(0,\frac7{64}\log p]$, then by \eqref{Eq:Chebyshevpoly} and \eqref{betadef}, for $\ell\leq L_q$, we have
\begin{align*}
X_\ell(\theta_k(q))&=\sum_{j=0}^\ell e^{(\ell-2j)\beta_k(q)}\\
&=\sum_{j=0}^{[\ell/2]} (e^{(\ell-2j)\beta_k(q)}+e^{-(\ell-2j)\beta_k(q)})-1+\ell-2[\ell/2]\\
&=X_\ell(0)+O(\ell^3|\beta_k(q)|^2),
\end{align*}
where $[\ell/2]$ means the integral part of $\ell/2$. Hence, by \eqref{al} and \eqref{bound for g}, we obtain that
$$
\mathfrak{g}_q(\theta_k(q))=\mathfrak{g}_q(0)+O\left(\sum_{\ell=0}^{L_q} |a_{\ell,q}|\ell^3|\beta_k(q)|^2\right)=1+O\left(L_q^4|\beta_k(q)|^2\right).
$$
Noting that $e^{in\pi}=e^{-in\pi}$ for any $n\in\mathbb{Z}$, by similar arguments, if $\theta_k(q)\in \pi+i(0,\frac7{64}\log p]$, then 
$$
\mathfrak{g}_q(\theta_k(q))=\mathfrak{g}_q(\pi)+O\left(L_q^4|\beta_k(q)|^2\right)\leq 1+O\left(L_q^4|\beta_k(q)|^2\right).
$$
Here we have used \eqref{bound for g} in the last step. 
On the other hand, by \eqref{eq:2} and the Taylor expansion, we have
\begin{equation*}
|\lambda_k(q^H)|
=\sum_{h=0}^{[H/2]}(e^{(H-2h)|\beta_k(q)|}+e^{-(H-2h)|\beta_k(q)|})-1+H-2[H/2]
\geq H+1+H^2|\beta_k(q)|^2.
\end{equation*}
Then there exists $T_2>0$ such that for any $T>T_2$ and $|\beta_k(q)|<\frac{L_q+2\log H}{2H-L_q-1}$, we have
$$
\mathfrak{g}_q(\theta_k(q))\leq 1+O\left(L_q^4|\beta_k(q)|^2\right)\leq (1+H|\beta_k(q)|^2/2)^2\leq \frac{\lambda_k(q^H)^2}{H^2}
$$
and so combining with \eqref{bound for g} and \eqref{gpositive}, we have 
$$
\mathcal{G}(u_k)
\le
\left(\mathfrak{g}_q(\theta_k(q))-\frac{\lambda_k(q^H)^2}{H^2}\right)\prod_{\substack{p\le z\\ p\ne q}} \mathfrak{g}_p(\theta_k(p))
\le 0.
$$

Let $\A$ be the set of $u_k\in \hk$ such that $\theta_k(p)\in[0,2\pi \delta_p)$ for all primes $p\le z$,
and denote by $\mathbf{1}_{\A}$ the indicator function of $\A$. For all $u_k\in\A$ we have $\mathcal{G}(u_k)\le 1$ since $0\leq \mathfrak{g}_p(\theta)\leq 1$ for all $\theta\in \mathbb{R}$ and $p\leq z$ by \eqref{bound for g}. Assume $T$ is sufficiently large such that $T>\max(T_1,T_2)$.
Combining with the above three cases, we have
\begin{equation}\label{Eq:Compare1G}
\sum_{u_k\in\mathcal{H}_T}  \alpha_k \mathbf{1}_{\A}(u_k)
\ge
\sum_{u_k\in\mathcal{H}_T}  \alpha_k \mathcal{G}(u_k).
\end{equation}

\section{Proof of Theorem \ref{Thm:Main}}

We follow the approach in \cite{La26} with a refinement modifying the sieve used in his work to address the exceptional eigenvalues (due to the absence of GRC). Writing $p_j$ as the $j$th prime number and putting $J=\pi(z)$, we obtain by \eqref{Eq:ChebyshevExpansion} and \eqref{Eq:DefG}
\begin{align*}
\mathcal{G}(u_k)
&=
\prod_{p\le z}
\left(\sum_{\ell=0}^{L_p} a_{\ell,p} X_\ell(\theta_k(p))\right)
-
\sum_{q\le z}\left(\varepsilon+\frac{\lambda_k(q^H)^2}{H^2}\right)
\prod_{\substack{p\le z\\ p\ne q}}
\left(\sum_{\ell=0}^{L_p} a_{\ell,p} X_\ell(\theta_k(p))\right)\\
&=
\sum_{0\le \ell_1\le L_{p_1}}\cdots\sum_{0\le \ell_J\le L_{p_J}}
\prod_{j=1}^{J}
a_{\ell_j,p_j}X_{\ell_j}(\theta_k(p_j))\\
&\qquad
-
\sum_{t=1}^{J}\left(\varepsilon+\frac{\lambda_k(p_t^H)^2}{H^2}\right)
\sum_{\substack{0\le \ell_j\le L_{p_j}\\ 1\le j\le J\\ j\ne t}}
\prod_{\substack{j=1\\ j\ne t}}^{J}
a_{\ell_j,p_j}X_{\ell_j}(\theta_k(p_j)).
\end{align*}
Note that for any prime $p$,
$$
\lambda_k(p^H)^2=\sum_{h=0}^{2H}b_h\lambda_k(p^h),
$$
where
\begin{equation}\label{bh}
b_h=\frac{2}{\pi}\int_0^\pi X_H^2(\theta)X_h(\theta)(\sin\theta)^2\,d\theta.
\end{equation}
Inserting the above formula, and using \eqref{CP} for all $\ell$ gives
\begin{align}\label{sumofG}
\frac{\pi^2}{T^2}\sum_{u_k\in\hk}  \alpha_k \mathcal{G}(u_k)
&=\frac{\pi^2}{T^2}\sum_{\substack{0\le \ell_j\le L_{p_j}\\ 1\le j\le J}}
\prod_{j=1}^{J}
a_{\ell_j,p_j}\sum_{u_k\in\hk}   \alpha_k \lambda_k\left(\prod_{j=1}^{J} p_j^{\ell_j}\right)
\nonumber\\
& \quad \quad \quad  -
\frac{\varepsilon\pi^2}{T^2}
\sum_{t=1}^{J}
\sum_{\substack{0\le \ell_j\le L_{p_j}\\ 1\le j\le J\\ j\ne t}}
\prod_{\substack{j=1\\ j\ne t}}^{J}
a_{\ell_j,p_j} \sum_{u_k\in\hk}   \alpha_k \lambda_k\bigg(\prod_{\substack{j=1\\ j\ne t}}^{J} p_j^{\ell_j}\bigg)
\nonumber\\
& \quad \quad \quad  -
\frac{\pi^2}{T^2H^2}
\sum_{t=1}^{J}\sum_{h=0}^{2H}b_h
\sum_{\substack{0\le \ell_j\le L_{p_j}\\ 1\le j\le J\\ j\ne t}}
\prod_{\substack{j=1\\ j\ne t}}^{J}
a_{\ell_j,p_j} \sum_{u_k\in\hk}  \alpha_k  \lambda_k\bigg(p_t^h\prod_{\substack{j=1\\ j\ne t}}^{J} p_j^{\ell_j}\bigg).
\end{align}
Then we apply Lemma \ref{lem3.1} and \eqref{al} to the right hand side of \eqref{sumofG} and obtain
\begin{align*}
\frac{\pi^2}{T^2}\sum_{u_k\in\hk}  \alpha_k \mathcal{G}(u_k)= &\prod_{j=1}^{J}a_{0,p_j}
+O\left(T^{-\kappa_0}\prod_{p\le z} p^{L_p\eta_0} (L_p+1)
\right)\nonumber\\
&+\varepsilon
\sum_{t=1}^{J}
\prod_{\substack{j=1\\ j\ne t}}^{J}a_{0,p_j}+O\left(\varepsilon T^{-\kappa_0}\sum_{t=1}^{J}
\prod_{\substack{j=1\\ j\ne t}}^{J}p^{L_p\eta_0} (L_p+1)
\right)\\
&+\frac{b_0}{H^2}\sum_{t=1}^{J}
\prod_{\substack{j=1\\ j\ne t}}^{J}a_{0,p_j}+O\left(\frac1{H^2T^{\kappa_0}}\sum_{t=1}^{J}\sum_{h=0}^{2H}|b_h| p_t^{h\eta_0}
\prod_{\substack{j=1\\ j\ne t}}^{J}p^{L_p\eta_0} (L_p+1)
\right),
\end{align*}
where $\kappa_0=1-10^{-2026}$ and $\eta_0=\frac14+10^{-2026}$. 
By \eqref{bh}, we find that
$$
b_0=1\qquad\textrm{and}\qquad b_h\ll h.
$$
Hence, we have
\begin{align}\label{Eq:EstimateAverageG}
&\frac{\pi^2}{T^2}\sum_{u_k\in\hk}  \alpha_k  \mathcal{G}(u_k)\nonumber\\
&= \prod_{j=1}^{J}a_{0,p_j}+O\left((\varepsilon +1/H^2)\sum_{t=1}^{J}
\prod_{\substack{j=1\\ j\ne t}}^{J}a_{0,p_j}+T^{-\kappa_0}Jz^{2H\eta_0}\prod_{p\le z} p^{L_p\eta_0} (L_p+1)
\right).
\end{align}
Now we choose
$$
z=H =\left[\frac{\log  T}{2026m(\log \log T)} \right]
$$
so it infers that $J \sim  \log T/2026m(\log\log T)^2$. Then we have (see \cite[Page 7]{La26} for the proof)
\begin{align*}
a_{0,p_j} \gg \frac{1}{(\log\log T)^6},\ \textrm{for all}\ 1\leq j\leq J.
\end{align*}
Using these estimates in \eqref{Eq:EstimateAverageG} and recalling that $\ep\asymp (\log T)^{-\pi/2}$ we obtain
$$
\frac{\pi^2}{T^2}\sum_{f\in\hk}   \alpha_k \mathcal{G}(u_k)
\sim \prod_{j=1}^{J}a_{0,p_j}
\gg
\exp\left(-\frac{ (\log T)\log\log\log T}{2m(\log\log T)^2}\right).
$$
Next we combine this with the inequality \eqref{Eq:Compare1G} and the Weyl law \eqref{weyl} and \eqref{alphak} to derive
\begin{equation*}\label{Eq:LowerBoundSetA}
|\mathcal{A}| \gg |\hk| \exp\left(-\frac{ (\log T)\log\log\log T}{m(\log\log T)^2}\right).
\end{equation*}

Now let $u_k\in \A$ and suppose that $p^a\leq z$ for some prime $p$ and positive integer $a$. Then $p\leq z$ and $a\leq (\log z)/\log p\leq (\log \log T)/\log p.$
On the other hand, by \eqref{defL}, we have
$$
\lambda_{\sym^m u_k}(p^a)=\sum_{j_0+j_1+\cdots+j_{m}=a}\alpha_{u_k,1}(p)^{a m-2\sum_{\ell=0}^m\ell j_\ell}=\sum_{h =0}^{am} c_h \cos h\theta_k(p),
$$
where $c_{am}=2$, $c_h\geq0$ for $h=0,2,\ldots,am$ and we have used the fact $\lambda_{\sym^m u_k}(p^a)$ is a symmetric polynomial in $\alpha_{u_k,1}(p)$ and $\alpha_{u_k,1}(p)^{-1}$.
Therefore, for $h=0,2,\ldots,am$ and
$0\leq \theta_k(p)\leq2\pi \delta_p$ (by our assumption that $u_k\in \A$), we have $$0\leq h \theta_k(p)\leq am \theta_k(p)\leq 4 m \pi\delta_p (\log\log T)/\log p\leq \frac{\pi}{4},$$
which implies that
$$
\lambda_{\sym^m u_k}(p^a)>0.
$$
This concludes the proof.

\end{document}